\documentclass[11pt, reqno]{amsart}
\usepackage{lmodern}

\usepackage{amsmath,amsfonts,amssymb,amsthm,booktabs,color,epsfig,graphicx,hyperref,url}

\usepackage[table]{xcolor}

\usepackage{natbib}
\usepackage{dsfont}
\usepackage{upgreek, mathrsfs}
\usepackage{tikz}

\usepackage{rotating}
\usepackage{array} 
\usepackage{longtable}
\usepackage{enumitem}
\usepackage{calc} 
\allowdisplaybreaks

\usepackage{setspace}
\newcommand{\bbz}{\mathbb{Z}}
\newcommand{\bbn}{\mathbb{N}}

\newcommand{\bfZ}{\mathbf{Z}}

\newcommand{\bbe}{\mathbb{E}}
\newcommand{\var}{\mathbb{V}}
\newcommand{\bbp}{\mathbb{P}}
\newcommand{\cov}{\text{Cov}}
\newcommand{\corr}{\text{Corr}}

\newcommand{\pgf}{\operatorname{pgf}}

\newcommand{\brackets}[1]{\left( #1 \right)}
\newcommand{\bbrackets}[1]{\big( #1 \big)}

\theoremstyle{plain}
\newtheorem{theorem}{Theorem}[section]

\theoremstyle{definition}

\newtheorem{remark}[theorem]{Remark}

\begin{document}
\title[]{INARMA Models for Count Random Fields --- a Survey} 
\author[A. Silbernagel]{Angelika Silbernagel\textsuperscript{\textasteriskcentered}}
\author[C.H. Wei\ss{}]{Christian H. Wei\ss{}\textsuperscript{\dag}}

\address{\textsuperscript{\textasteriskcentered} Department of Mathematics and Statistics, Helmut Schmidt University, Hamburg, Germany. \newline
ORCID: \href{https://orcid.org/0009-0002-6993-244X}{\nolinkurl{0009-0002-6993-244X}}.
}
\email{silbernagel@hsu-hh.de}

\address{\textsuperscript{\dag} Department of Mathematics and Statistics, Helmut Schmidt University, Hamburg, Germany. \newline
ORCID: \href{https://orcid.org/0000-0001-8739-6631}{\nolinkurl{0000-0001-8739-6631}}.
}
\email{weissc@hsu-hh.de}


\begin{abstract}
The thinning-based integer-valued autoregressive moving-average (INARMA) models are popular for count time series. Recently, types of INARMA models have also been developed for count random fields, i.e., for spatial count data located on a regular two-dimensional grid. This article provides a comprehensive survey on existing INARMA random fields, covering approaches with different thinning operators, first- and higher-order models, as well as unilateral and multilateral model structures.
\end{abstract}
\maketitle

\section{Introduction}
\label{section: Introduction}

Spatial data located on a regular two-dimensional grid (also known as lattice data or raster data) are of interest in many fields of practice, covering agriculture \citep{tabandeh24}, biology \citep{chutoo21}, epidemiology \citep{yang26}, and political sciences \citep{adaemmer26}. The stochastic process generating such grid data is commonly referred to as a random field (also spatial process in the plane). Real-world grid data (like in the aforementioned applications) are typically characterized by spatial dependencies that have to be accounted for when selecting a stochastic model for the underlying random field. With a few exceptions, the research first focused on real-valued and continuously distributed random fields. In particular, ``plane counterparts'' to the classical autoregressive (AR) and moving-average (MA) models known from time series analysis have been developed and extensively discussed several decades ago, see, e.g., \citet{basu93, besag74, haining78, whittle54}. 

\begin{figure}[t]
    \centering
    (a)\hspace{-3ex}\includegraphics[viewport=30 45 185 245, clip=, scale=0.65]{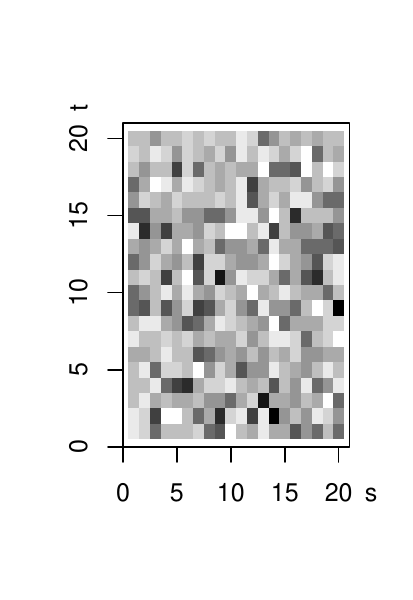}
    \qquad
    (b)\includegraphics[viewport=30 55 360 255, clip=, scale=0.65]{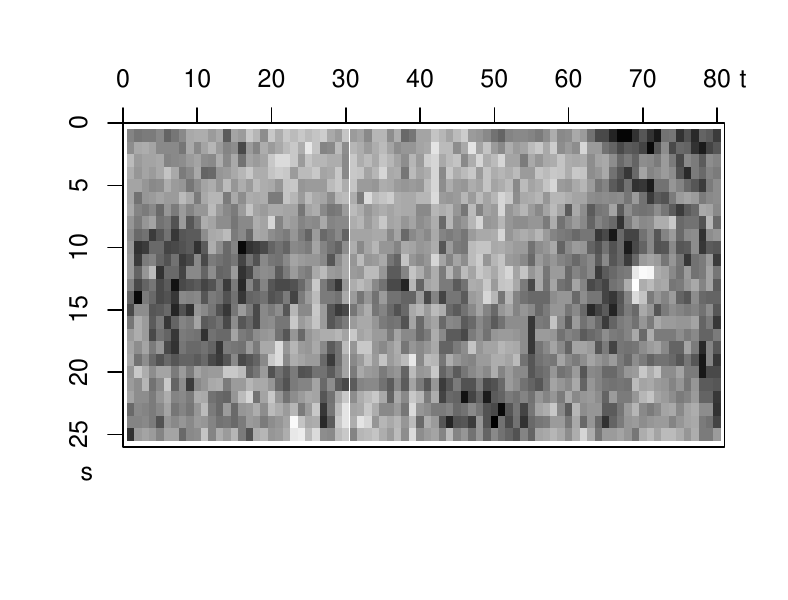}
    \caption{Illustrative data examples: (a) yeast count data of \citet{student06}, and (b) wheat yields data of \citet{iyer42}, where increasing counts are represented by darker shades of gray.}
\label{fig: examples}
\end{figure}

\smallskip
However, the grid data observed in real-world applications can sometimes not be explained by continuously distributed random fields, like it happens in the aforementioned data examples: \citet{tabandeh24} analyze counts of striga infestation in a millet field, \citet{chutoo21} consider the numbers of trees of the species Beilschmiedia pendula in a tropical rain forest, \citet{yang26} model counts of COVID-19 infections in China, and \citet{adaemmer26} investigate the numbers of war-related fires in Ukraine. Also the data examples being plotted in Fig.~\ref{fig: examples} (which will be referred to again later in this article) are cases of count grid data: part~(a) shows counts of yeast cells over 1~mm\textsuperscript{2} divided into $20\times 20$ squares \citep[p.~355]{student06}, while~(b) shows counts of half-ounce units concerning the yield of wheat in an experimental field with $25\times 80$ plots \citep{iyer42}. 
So in all these examples, the outcomes are discrete-valued counts as generated by count random fields $(X_{s,t}) = (X_{s,t})_{s,t\in\bbz=\{\ldots, -1,0,1, \ldots\}}$, i.e., random fields having the non-negative integers $\bbn_0=\{0,1,\ldots\}$ as their range.

\smallskip
The development of stochastic models for discrete-valued random fields has been widely ignored for a long time. The few early proposals include regression-type approaches such as the so-called ``auto-Poisson model'' by \citet{besag74,kaiser97} on the one hand, and the ``Poisson/gamma random field model'' of \citet{wolpert98} on the other hand. The first belongs to the class of observation-driven models, the latter to the parameter-driven ones, with both classes having their pros and cons \citep{koopman16}. For example, observation-driven models typically lead to closed-form expressions for the likelihood function and thus allow for maximum likelihood estimation by direct numerical optimization. However, there are also computationally feasible estimation approaches for parameter-driven models, such as the integrated nested Laplace approximation (INLA) approach being applicable to latent Gaussian random-field models \citep{martino20}. Pros and cons are also noted for various Gaussian random field models, which allow for a tailor-made choice of count marginal distributions on the one hand, but which suffer from a challenging parameter estimation on the other hand, see \citet{masarotto12,kazianka13,hughes15,han20} and the references therein for a discussion.

\smallskip
Discrete-valued counterparts to another observation-driven approach for random fields, the popular ARMA models mentioned before, have not been discussed in the research literature for a long time. In fact, the first integer-valued AR-type (INAR) model for count random fields has been introduced only about~$15$ years ago, and integer-valued MA-type (INMA) random fields have been considered even more recently. These INARMA models for count random fields have been developed as ``planar counterparts'' of classical INARMA models known from time series analysis, see \citet{weiss18} for a comprehensive survey on the latter. Like their time-series counterparts, INARMA models for random fields have a well-interpretable data-generating mechanism and a linear dependence structure that is described by kinds of Yule--Walker (YW) equations. Both of these aspects make the models potentially attractive to users. However, due to their relatively short history, the INARMA random fields are not yet widely known among practitioners. For this reason, and also to identify open questions for future research, the present paper aims to provide a comprehensive overview of the existing contributions on INARMA random fields.
In accordance to the historical development, Sect.~\ref{section: The Origins: The First-order INAR Random Field} starts with a presentation of the first INARMA model for count random fields, namely the first-order INAR model proposed by \citet{ghodsi12}. Here and in the subsequent sections, we complement the model definition with fundamental stochastic properties of the respective models, such as marginal mean and variance as well as the spatial autocorrelation structure. Sect.~\ref{section: Modifications Based on Alternative Thinning Concepts and Mixtures} deals with modifications of this first INAR model in terms of alternative thinning concepts and mixtures applied to the recursion formula. Extensions to multilateral structure and higher model orders are discussed in Sects.~\ref{section: INAR Extensions to Multilateral Structure} and~\ref{section: INAR Extensions to Higher Orders}, while Sect.~\ref{section: INMA Random Fields} considers the INMA model for count random fields. 
A short conclusion and an outlook on future research opportunities is presented in Sect.~\ref{section: Conclusions and Outlook}.

\section{The Origins: The First-order INAR Random Field}
\label{section: The Origins: The First-order INAR Random Field}

The first INAR model for count random fields has been proposed by \citet{ghodsi12}, and it was referred to as the ``first-order spatial INAR (SINAR$(1,1)$)'' model by the authors. It does not only serve as a planar analog of the INAR$(1)$ time-series model, but also as an integer-valued counterpart to the continuously distributed spatial AR model of \citet{basu93}, the latter being a spatial analog of the classical AR time-series model. Since there are quite different types of spatial data (and corresponding models) in the literature, \citet{sil_wei_26} later referred to the SINAR$(1,1)$ process as the INAR$(1,1)$ \emph{random field} instead, in order to indicate the grid structure of the generated outcomes. In what follows, to prevent any confusion, we uniquely adapt the latter terminology. 

\smallskip
In order to properly introduce the INAR$(1,1)$ random field, first we need to discuss the probabilistic operator of binomial thinning, which goes back to \citet{steutel79} and constitutes a fundamental tool in the modeling of integer-valued time series. Let $X$ be a discrete random variable with range~$\bbn_0$, and let $\alpha \in (0,1)$ be a constant. We say that the random variable $\alpha \circ X := \sum^X_{i=1} Z_i$ arises from $X$ by \emph{binomial thinning}, where the \emph{counting series} $(Z_i)_{i\in\bbn=\{1,2,\ldots\}}$ is a sequence of independent and identically distributed (i.i.d.)\ Bernoulli random variables with $\bbp(Z_i=1)=\alpha$, which is independent of $X$. So by construction, the count random variable $X$ is ``reduced'' to $\alpha \circ X$ as the latter can only take values between $0$ and $X$. Note that given $X$, $\alpha \circ X$ follows a binomial distribution, that is, $\alpha \circ X \,|\, X \sim \text{Bin}(X, \alpha)$. Considering a population of size $X$, the thinned random variable $\alpha \circ X$ can be understood as, e.g., the number of survivors, where each individual $Z_i$ survives independently of each other and with probability $\alpha$, see \citet{weiss18}. The boundary values $\alpha\in\{0,1\}$ can be included via the conventions $0 \circ X := 0$ and $1 \circ X:=X$. Now, the \emph{INAR$(1,1)$ random field} is defined by the recursion
\begin{equation}
    \label{eq: model}
    X_{s,t} = \alpha_{10} \circ X_{s-1,t} + \alpha_{01} \circ X_{s,t-1} + \alpha_{11} \circ X_{s-1,t-1} + \varepsilon_{s,t}\,,
\end{equation}
with the dependence parameters satisfying
\begin{equation}
    \label{eq: parameter condition}
    \alpha_{10}, \alpha_{01}, \alpha_{11} \in [0,1)
    \quad\text{and}\quad
    \alpha_{10}+\alpha_{01}+\alpha_{11}<1\,.
\end{equation}
The count innovations $(\varepsilon_{s,t})=(\varepsilon_{s,t})_{s,t\in\bbz}$ are assumed to be i.i.d.\ with mean $\mu_\varepsilon$ and variance $\sigma^2_\varepsilon$, and the thinning operators at point $(s,t)$ appearing in \eqref{eq: model} are assumed to be performed independently of each other and of $(\varepsilon_{i,j})$. Moreover, if we define $\mathcal{P}\!_{s,t}=\{X_{s-k,t-l} : k\geq1 \text{ or } l\geq 1\}$ for the ``past'' corresponding to ``time'' $(s,t)$, then the innovation $\varepsilon_{s,t}$ at ``time'' $(s,t)$ is assumed to be independent of $\mathcal{P}\!_{s,t}$ for all $s,t\in\bbz$. Finally, in order to being able to derive the autocorrelation function (ACF) $\rho(k,l) := \corr(X_{s,t}, X_{s+k,t+l})$, we also have to assume that all thinnings performed at ``time'' $(s,t)$ are independent of $\mathcal{P}\!_{s,t}$, see \citet{sil_wei_26} for a discussion. 

\smallskip
The model in \eqref{eq: model} is proven to be ergodic in \citet{ghodsi15}. Under stationarity assumptions, its mean equals $\mu_X = \bbe(X_{s,t}) = \mu_\varepsilon/(1-\alpha_{10}-\alpha_{01}-\alpha_{11})$, see \cite{ghodsi12}, and its autocovariance function (ACvF) $\gamma(k,l):= \cov(X_{s,t}, X_{s+k,t+l})$ satisfies the YW equations
\begin{equation} 
    \label{eq: YW equations}
    \begin{split}
    \gamma(k,l) &= \alpha_{10}\, \gamma(k-1,l) + \alpha_{01}\, \gamma(k,l-1) + \alpha_{11}\, \gamma(k-1,l-1) \
    \text{ for } k \geq 1 \text{ or } l \geq 1\,, 
    \\[1ex]
    \gamma(k,l) &= \alpha_{10}\, \gamma(k+1,l) + \alpha_{01}\, \gamma(k,l+1) + \alpha_{11}\, \gamma(k+1,l+1) \
    \text{ for } k \leq -1 \text{ or } l \leq -1\,,
    \end{split}
\end{equation} 
see \citet{ghodsi12,sil_wei_26}. If, moreover, at most one autoregressive coefficient $\alpha_{ij}$, $(i,j)\in\{(1,0),(0,1),(1,1)\}$, is zero, then the ACvF $\gamma(k,l)$ in the region ($k\geq 0$ and $l \leq 0$) has the closed-form expression
\begin{equation}
\label{eq: sacf main}
    \gamma(k,l) = \gamma(0,0)\, \lambda^{k}\, \eta^{-l}\,,
\end{equation}
with $\lambda, \eta \in (0,1)$ given by
\begin{equation}
    \label{eq: lambda eta}
    \begin{split}
        &\lambda = \frac{(1+\alpha_{10}^2-\alpha_{01}^2-\alpha_{11}^2)-\sqrt{(1+\alpha_{10}^2-\alpha_{01}^2-\alpha_{11}^2)^2-4(\alpha_{10}+\alpha_{01}\alpha_{11})^2}}{2(\alpha_{10}+\alpha_{01}\alpha_{11})} \\
        &\quad \text{ and } \quad
        \eta=\frac{\alpha_{01}+\alpha_{11}\lambda}{1-\alpha_{10}\lambda}\,.
    \end{split}
\end{equation}
From this, it can be shown that the variance satisfies
\begin{equation}
    \label{eq: gamma00}
        \sigma_X^2 = \gamma(0,0) = \frac{\mu_X\bbrackets{\alpha_{10}(1-\alpha_{10})  + \alpha_{01}(1-\alpha_{01}) + \alpha_{11}(1-\alpha_{11})} + \sigma_\varepsilon^2}{\sqrt{(1+\alpha_{10}^2-\alpha_{01}^2-\alpha_{11}^2)^2-4(\alpha_{10}+\alpha_{01}\alpha_{11})^2}}\,.
\end{equation}

\citet{ghodsi12} consider the YW method for parameter estimation, while \citet{ghodsi15} discusses the conditional maximum likelihood (CML) method in the specific case of the Poisson INAR$(1,1)$ model (i.e., the innovations~$\varepsilon_{s,t}$ are assumed to follow a Poisson distribution) and derive the asymptotic distribution of the estimator. In both cases, the estimation methods were used to fit the INAR$(1,1)$ model to the yeast count data being plotted in Figure~\ref{fig: examples}\,(a). The conditional least squares (CLS) estimation method is considered by \citet{sassi23}, where the authors show strong consistency and asymptotical normality (without distributional assumption on the innovations). We conclude by noting that despite having Poisson innovations~$\varepsilon_{s,t}$, the Poisson INAR$(1,1)$ model does not lead to Poisson observations~$X_{s,t}$ (in fact, the distribution of~$X_{s,t}$ is not known). Poisson observations, by contrast, can be achieved with the combined INAR (CINAR) model of \citet{wei_sil_26}, which is presented in Sect.~\ref{section: INAR Extensions to Higher Orders} below.

\section{Modifications Based on Alternative Thinning Concepts and Mixtures}
\label{section: Modifications Based on Alternative Thinning Concepts and Mixtures}

The INAR$(1,1)$ random field model of \citet{ghodsi12} has been modified and extended in many subsequent works.
\citet{ghodsi24}, \citet{tabandeh24} and \citet{yang26} equip the model with different types of thinning operators. 
\citet{ghodsi24} opt for the \emph{generalized binomial thinning operator} denoted by ``$\circ_\theta$'' with $\theta\in[0,1]$, which is a mixture of two binomial distributions:
\[
    \alpha \circ_\theta X \,|\, X\ \sim\ (1-\alpha)\, \text{Bin}\big(X,\alpha(1-\theta)\big) + \alpha\, \text{Bin}\big(X, \alpha(1-\theta)+\theta\big)\,.
\]
\citet{tabandeh24}, by contrast, use a \emph{random coefficient thinning operator} in the model recursion \eqref{eq: model}, i.e., the fixed parameters $\alpha_{10}, \alpha_{01}, \alpha_{11}$ are replaced by random parameters $\alpha_{10}{}^{(s,t)}, \alpha_{01}{}^{(s,t)}, \alpha_{11}{}^{(s,t)}$ with means $\alpha_{10}, \alpha_{01}, \alpha_{11}$, respectively. The resulting sequences $\alpha_{ij}{}^{(s,t)}$ are assumed to be i.i.d.\ on $[0,1)$ for fixed $(i,j)$. Finally, \citet{yang26} consider the \emph{negative-binomial (NB) thinning operator} ``$\ast$'', which is defined similarly to the binomial thinning operator but uses a geometrically distributed counting series: $\alpha \ast X := \sum^X_{i=1} Z_i$ with $Z_i\sim\textup{Geom}\big(1/(1+\alpha)\big)$. This way, the counting series is allowed to have the full range $\bbn_0$ instead of just $\{0,1\}$ as it is in the case for the binomial thinning operator. Thus, $\alpha \ast X$ may become larger than $X$, so the interpretation of surviving individuals from a population of size $X$ is no longer applicable. Instead, it can be understood as a reproduction mechanism, where $Z_i$ describes the number of children generated by the $i$-th individual of the population (see \citealp[Sect.~3.2]{weiss18}). The operator $\alpha \ast X$ is conditionally NB-distributed, $\alpha \ast X \,|\, X \sim \text{NB}(X, 1/(1+\alpha))$, and useful for modeling count data with overdispersion (see \citealp{yang26}). Apart from changing the thinning operator, \citet{yang26} specify the innovations $\varepsilon_{s,t}$ to be geometrically distributed, and show the existence of a stationary and ergodic random field satisfying their model definition. Interestingly, it seems that the case of Poisson thinning (together with Poisson innovations) has not been discussed so far. As argued by \citet{weiss15}, this combination would offer a bridge to another popular class of models from count time series analysis, namely the integer-valued generalized AR conditional heteroscedasticity (INGARCH) models, and is, thus, recommended for future research, see Section~\ref{section: Conclusions and Outlook}.

\smallskip
The mean, variance, and ACF of the respective stationary INAR modifications are derived in the aforementioned papers. Interestingly, mean and ACF agree with those of \eqref{eq: model} in all three cases, whereas the variances of the models by \citet{ghodsi24} and \citet{tabandeh24} deviate from \eqref{eq: gamma00} due to the different thinning concepts. 
This is caused by the fact that, in analogy to the conditionally linear AR processes studied by \citet{grunwald00}, the conditional mean satisfies the linear structure (in addition to stationarity)
\begin{equation*}
    \bbe(X_{s,t} \,|\, \mathcal{P}_{s,t}) 
    = \alpha_{10}\, X_{s-1,t} + \alpha_{01}\, X_{s,t-1} + \alpha_{11}\, X_{s-1,t-1} + \mu_\varepsilon.
\end{equation*}
Then, it follows by repeated application of the tower property that the (non-conditional) mean and YW equations are given by $\mu_X = \mu_\varepsilon/(1-\alpha_{10}-\alpha_{01}-\alpha_{11})$ and \eqref{eq: YW equations}, respectively. The variance 
\begin{align*}
    \var(X_{s,t}) &=  \var(\bbe(X_{s,t} \,|\, \mathcal{P}_{s,t})) + \bbe(\var(X_{s,t} \,|\, \mathcal{P}_{s,t})) \\
    &= (\alpha_{10}^2 + \alpha_{01}^2 + \alpha_{11}^2)\, \sigma_X^2 + 2\alpha_{10} \alpha_{01}\, \gamma(1,-1) + 2 \alpha_{10}\alpha_{11}\, \gamma(0,1) \\
    &\qquad \qquad + 2\alpha_{01}\alpha_{11}\, \gamma(1,0) + \bbe(\var(X_{s,t} \,|\, \mathcal{P}_{s,t})),
\end{align*}
however, also depends on the conditional variance, which differs between the different INAR random fields. Finally, 
in all three papers, YW, CLS, and CML methods are considered for parameter estimation. Additionally, \citet{yang26} also use the weighted CLS method, and derive the asymptotic properties of all estimators under consideration. 

\medskip
Another extension of the INAR random field model has been proposed by \citet{ghodsi25}, which combines the so-called Pegram operator and the binomial thinnings within \eqref{eq: model}. The \emph{Pegram operator} ``$\star$'' mixes two random variables $X$ and $Y$ with proportions $\phi$ and $1-\phi$, and it is expressed as $Z=(\phi, X) \star (1-\phi, Y)$. 
Note that \citet{ghodsi25} denote the Pegram operator by ``$\ast$'', but we use the symbol ``$\star$'' to distinguish it from the NB-thinning operator.
Then, the model recursion of \citet{ghodsi25} is given by
\[
    X_{s,t} = (\phi,\, \alpha_{10} \circ X_{s-1,t} + \alpha_{01} \circ X_{s,t-1} + \alpha_{11} \circ X_{s-1,t-1}) \star (1-\phi,\, \varepsilon_{s,t})\,,
\]
which is a mixture of the AR-part and the innovation~$\varepsilon_{s,t}$. \citet{ghodsi25} state that the model is capable of generating multimodal data. They derive its fundamental properties and discuss the CML method for parameter estimation. In particular, the ACF agrees with \eqref{eq: sacf main} if $\alpha_{ij}$ is replaced by $\phi\, \alpha_{ij}$, while the formulas for mean and variance slightly differ from those in Sect.~\ref{section: The Origins: The First-order INAR Random Field}. At this point, it should be noted that also the aforementioned CINAR model of \citet{wei_sil_26} is a kind of mixture model, see Sect.~\ref{section: INAR Extensions to Higher Orders} for details.

\smallskip
Finally, \citet{chutoo21} extend the INAR$(1,1)$ random field from \eqref{eq: model} to the non-stationary case by allowing for covariate information to influence the mean of the innovation process. The authors note that closed-form expressions for the moments are difficult to obtain if the model is non-stationary. Therefore, they consider the CML method for parameter estimation as it does not require for moment expressions.

\section{INAR Extensions to Multilateral Structure}
\label{section: INAR Extensions to Multilateral Structure}

\begin{figure}[t]
    \centering
    \includegraphics[viewport=10 25 385 310, clip=, scale=0.6]{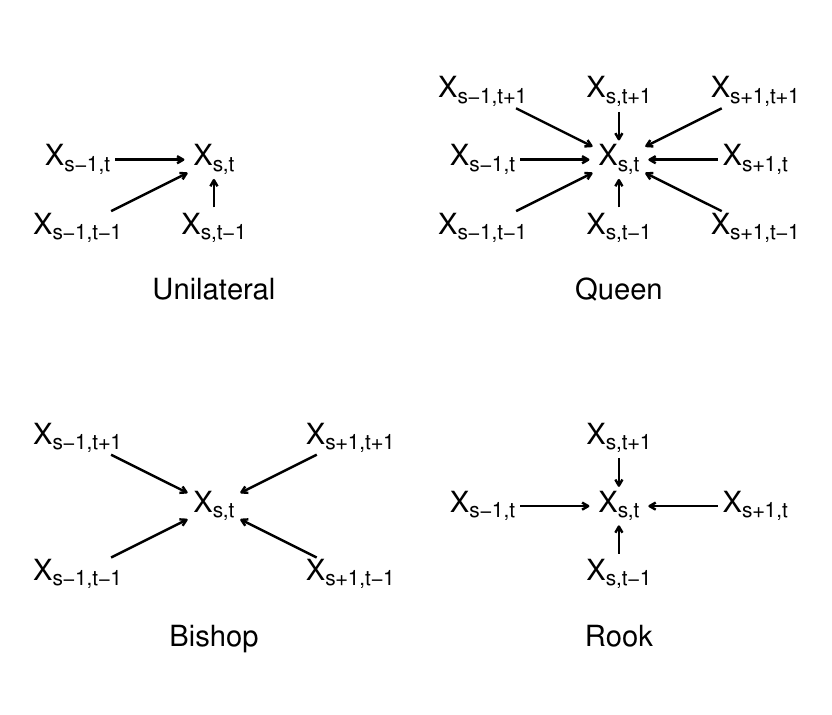}
    \caption{Different dependence structures for first-order models. The arrows indicate which observations influence the observation $X_{s,t}$. The Bishop, Rook, and unilateral models can be understood as special cases of the Queen model. Graphs adapted from \citet[Fig.~1]{karlis24}.}
    \label{fig: multilateral structure}
\end{figure}

\citet{karlis24} extend the unilateral structure of \eqref{eq: model}, where $X_{s,t}$ is explained solely by the ``past'' $X_{i,j}$ with $i\leq s$ \emph{and} $j\leq t$ (note the difference to $\mathcal{P}\!_{s,t}$), by different multilateral structures which allow for including ``future'' observations $X_{i,j}$ with $i\geq s$ or $j\geq t$. The different structures (for first-order models) are illustrated by Fig.~\ref{fig: multilateral structure}. As the Bishop and Rook structure are special cases of the Queen architecture (by setting the respective dependence parameters to zero), we discuss only the latter. Defining $\mathcal{S}^\ast:= \{(i,j) : i,j\in \{-1,0,1\}, (i,j) \neq (0,0)\}$ as the set of neighbors of $(0,0)$, the \emph{multilateral INAR$(1,1,1,1)$ random field} is defined by
\begin{equation}
    \label{eq: multilateral model}
    X_{s,t} = \sum_{(i,j)\in\mathcal{S}^\ast} \alpha_{ij} \circ X_{s-i,t-j}\ + \varepsilon_{s,t}
\end{equation}
with dependence parameters satisfying $\alpha_{ij} \in [0,1)$ and $\alpha_\bullet:=\sum_{(i,j)\in\mathcal{S}^\ast} \alpha_{ij} < 1$, where the independence assumptions are transferred from the unilateral model \eqref{eq: model}. Note that we recover the unilateral INAR$(1,1)$ model by setting $\alpha_{-1,-1},\alpha_{-1,0}, \alpha_{-1,1}, \alpha_{0,-1}, \alpha_{1,-1}$ equal to zero.  
\citet{karlis24} show that under stationarity assumptions, the marginal mean satisfies $\mu_X = \mu_\varepsilon/ (1-\alpha_\bullet)$. Expressions for variance and ACvF become rather sophisticated and are, thus, omitted here. We only cite the YW equations according to (C1) in \citet{karlis24}, 
\begin{equation*}
    \gamma(k,l) = \sum_{(i,j)\in\mathcal{S}^\ast} \alpha_{ij}\, \gamma(k-i,l-j) \qquad \text{for } k>0 \text{ or } l>0\,,
\end{equation*}
the derivation of which seems to require additional independence assumptions on the thinnings in analogy to \citet{sil_wei_26}, i.e., that the thinnings executed at ``time'' $(s,t)$ are independent of all other observations. 
\citet{karlis24} propose a CML estimator for the model parameters and prove its asymptotic normality by extending the asymptotic properties in \citet{ghodsi15}. Finally, let us briefly refer to two other multilateral extensions: in \citet{wei_sil_26}, it is shown how the recursion of the CINAR model could be extended to the multilateral case, while \citet{sil_wei_26b} develop a multilateral INMA model. These approaches shall be briefly discussed in Sects.~\ref{section: INAR Extensions to Higher Orders} and~\ref{section: INMA Random Fields}, respectively. 

\begin{remark}
\label{remark:multilateral_open_questions}
As also noted by the reviewer, the current research on the multilateral INAR model has not yet answered all the key questions. For example, \citet{karlis24} do not provide a formal existence proof for the random field defined by \eqref{eq: multilateral model}. Also its simulation appears to be demanding, where \citet{karlis24} refer to the Gibbs sampler approach described in Section~4.1 of \citet{jackson08}. 
There are also open questions concerning model fitting. For example, it is not clear if the likelihood function can be computed in closed form, or if pseudo-likelihood approximations are necessary. In addition, conditions for the identifiability of the model parameters appear to be necessary. As pointed out by the reviewer, the (degenerate) multilateral INAR models defined by $X_{s,t}=\alpha\circ X_{s-1,t}+\epsilon_t$ and $X_{s,t}=\alpha\circ X_{s+1,t}+\epsilon_t$ with Poisson innovations are equivalent, because they correspond to the Poisson INAR$(1)$ time-series model that is known to be time-reversible \citep{schweer15}.
\end{remark}

\section{INAR Extensions to Higher Orders}
\label{section: INAR Extensions to Higher Orders}

All models discussed so far are first-order models. \citet{sil_wei_26} generalize the INAR$(1,1)$ model \eqref{eq: model} to higher-order autoregressions by considering the recursion
\begin{equation}
    \label{eq: generalized model}
    X_{s,t} = \sum_{(i,j)\in\mathcal{S}} \alpha_{ij} \circ X_{s-i,t-j} + \varepsilon_{s,t}
\end{equation}
with AR coefficients $\alpha_{ij}\in[0,1)$ for all $(i,j)\in\mathcal{S}$, where $\mathcal{S}:=\big\{(i,j) \,|\, 0 \leq i \leq p_1,\, 0 \leq j \leq p_2,\, (i,j) \neq (0,0)\big\}$ and $p_1,p_2\in\bbn$. Note that \eqref{eq: model} is recovered for $p_1=p_2=1$. The independence assumptions from the INAR$(1,1)$ random field model are retained. This generalized model is referred to as the \emph{INAR$(p_1,p_2)$ random field}. For convenience, we define again $\alpha_\bullet := \sum_{(i,j)\in\mathcal{S}} \alpha_{ij}$ as the sum of all AR coefficients, where $\alpha_\bullet<1$ is necessary for stationarity. Then, \citet{sil_wei_26} derive the mean and variance as $\mu_X = {\mu_\varepsilon}/{(1-\alpha_\bullet)}$ and
\begin{align*}
    \brackets{1-\sum_{(i,j)\in\mathcal{S}} \alpha_{ij}^2}\, \sigma^2_X\ &=\  \sum_{\substack{(i_1,j_1),(i_2,j_2) \in \mathcal{S}\\(i_1,j_1)\neq(i_2,j_2)}} \alpha_{i_1j_1} \alpha_{i_2j_2}\, \gamma(i_1-i_2,j_1-j_2) \\
    &\qquad \qquad  \ +\ \brackets{\sum_{(i,j)\in\mathcal{S}} \alpha_{ij}(1-\alpha_{ij})} \mu_X + \sigma^2_\varepsilon\,.
\end{align*}
Furthermore, the YW equations are given by
\begin{align}
\label{eq: sacf generalized recursion}
\begin{split}
    \rho(k,l) 
    &= \sum_{(i,j)\in\mathcal{S}} \alpha_{ij}\, \rho(k-i,l-j)
    \quad \text{for } k\geq1 \text{ or } l\geq 1\, , \\
    \rho(k,l) 
    &= \sum_{(i,j)\in\mathcal{S}} \alpha_{ij}\, \rho(k+i,l+j)
    \quad \text{for } k\leq-1 \text{ or } l\leq -1\, . 
\end{split}
\end{align}
INAR$(p_1,p_2)$ random fields have the drawback that closed-form expressions for the stationary marginal distribution are difficult to obtain, since multiple thinnings are involved at each point $(s,t)$, even if only the first-order case $p_1=p_2=1$ is considered. Therefore, \citet{wei_sil_26} propose a possible solution by adapting an approach developed by \citet{zhu06,weiss08} for time-series data to random fields. This solution combines the INAR$(1)$ \emph{time-series} model, which involves only one thinning per time, with a random-selection mechanism. Referring to this model as the CINAR one, the \emph{CINAR$(p_1,p_2)$ random field} is defined by
\begin{equation}
    \label{eq: CINAR recursion}
    X_{s,t} = \sum_{(i,j)\in\mathcal{S}} D_{s,t;i,j} \times (\alpha \circ X_{s-i,t-j}) + \varepsilon_{s,t}\, ,
\end{equation}
where $(\mathbf{D}_{s,t})_{s,t\in\bbz}$ is an i.i.d.\ random field of ``decision'' random variables $\mathbf{D}_{s,t}=(D_{s,t;i,j})_{(i,j)\in\mathcal{S}}\sim\text{MULT}(1;\phi_{01}, \dots, \phi_{p_1p_2})$ independent of the innovations $(\varepsilon_{s,t})$, where exactly one of the $|\mathcal{S}|$ components equals $1$ and all the others equal $0$. Moreover, the thinnings at ``time'' $(s,t)$ are performed independently of each other, of $(\varepsilon_{s,t})$, $(\mathbf{D}_{s,t})$, and the ``past'' $\mathcal{P}\!_{s,t}$, and $\varepsilon_{s,t}$ and $\mathbf{D}_{s,t}$ are independent of the ``past'' $\mathcal{P}\!_{s,t}$ as well. In a nutshell, this means that $X_{s,t}$ is equal to $\alpha \circ X_{s-i,t-j}+\varepsilon_{s,t}$ with probability $\phi_{i,j}$.

\smallskip
\citet{wei_sil_26} show that the CINAR random field is ergodic. Its probability generating function (pgf) $\pgf_X(u) = \bbe(u^{X_{s,t}})$, mean and variance coincide with those of the classical INAR$(1)$ time-series model. This implies that the CINAR random field \eqref{eq: CINAR recursion} can be equipped with Poisson and NB-marginal distributions (among others), which differs from the previous INAR random fields. 
The YW equations agree with \eqref{eq: sacf generalized recursion} if $\alpha_{ij}$ is replaced by $\alpha\, \phi_{ij}$.
Moreover, \citet{wei_sil_26} consider parameter estimation by means of YW, CLS, and CML methods, which were used to fit first- and second-order CINAR models to the wheat yields data from Figure~\ref{fig: examples}\,(b). 

\smallskip
In addition, \citet{wei_sil_26} briefly discuss possible future extensions of the basic CINAR model \eqref{eq: CINAR recursion}. 
First, they show that a \emph{multilateral} recursion of order $(p_1,p_2,q_1,q_2)\in\bbn^4$ could be defined by
    \begin{equation}
        \label{eq: multilateral CINAR recursion}
        X_{s,t}
        = \sum_{(i,j)\in\mathcal{S}^*} D_{s,t;i,j} \times (\alpha \circ X_{s-i,t-j}) + \varepsilon_{s,t}\, ,
    \end{equation}
where $\mathcal{S}^*:=\big\{(i,j) \,|\, -q_1 \leq i \leq p_1,\, -q_2 \leq j \leq p_2,\, (i,j) \neq (0,0)\big\}$, and where the ``past'' $\mathcal{P}\!_{s,t}$ is replaced by $\mathcal{A}_{s,t} = \{X_{u,v}: (u,v)\not=(s,t)\}$. As the advantages of the unilateral CINAR model (such as closed-form marginal distributions) would be preserved for this recursion (see \citealp{wei_sil_26}), it could become a competitor to the multilateral INAR model of \citet{karlis24}, recall Sect.~\ref{section: INAR Extensions to Multilateral Structure}. However, a detailed investigation of a multilateral CINAR structure is still pending (also recall Remark~\ref{remark:multilateral_open_questions}) and should be addressed in future research. The same applies to the second possible extension proposed by \citet{wei_sil_26}, where negative AR terms are enabled in \eqref{eq: CINAR recursion} by using a \emph{Tobit approach}, also see Sect.~\ref{section: Conclusions and Outlook} below.

\section{INMA Random Fields}
\label{section: INMA Random Fields}

Interestingly, the scientific literature largely focused on AR-type models for count random fields so far. For simulating spatially dependent count random fields, however, an MA-type model would be more convenient, as burn-in periods for INAR random fields are computationally expensive, especially for multilateral INAR random fields as considered by \citet{karlis24}. MA-type models have been considered in a few applications within simulation studies, e.g., in \citet{wei_kim_24}, but \citet{sil_wei_26b} were the first to introduce and study such a model from a theoretical perspective. They define the (unilateral) \emph{INMA random field} of order $(q_1,q_2)\in\bbn_0^2$ with $q_1+q_2\geq 1$ by the recursion
\begin{equation}
    \label{eq: model inma}
    X_{s,t} = \sum^{q_1}_{i=0} \sum^{q_2}_{j=0} \beta_{ij}\, \circ_{s,t}\, \varepsilon_{s-i,t-j}\,,
\end{equation}
where $\beta_{ij}\in [0,1]$ denote the model parameters; abbreviate $\beta_\bullet:=\sum^{q_1}_{i=0} \sum^{q_2}_{j=0} \beta_{ij}$. $(\varepsilon_{s,t})=(\varepsilon_{s,t})_{s,t\in\bbz}$ is again a sequence of i.i.d.\ count innovations with mean $\mu_\varepsilon$ and variance $\sigma^2_\varepsilon$, and all thinnings concerning \emph{different}~$\varepsilon_{s,t}$ are performed independently of each other and of~$\varepsilon_{s,t}$. Here, we have added an index ``$s,t$'' below the operator ``$\circ$'' in order to indicate that the thinning is performed at point $(s,t)$. 

\smallskip
\citet{sil_wei_26b} derive the INMA$(q_1,q_2)$'s marginal distribution as well as the spatial dependence structure under stationarity assumptions. In particular, mean, variance, and pgf of $X_{s,t}$ are given by 
\begin{equation}
\label{eq: INMA mean var pgf}
\begin{split}
    \mu_X &= \mu_\varepsilon\, \beta_\bullet\, ,
    \quad
    \sigma_X^2 = \mu_X + (\sigma^2_\varepsilon - \mu_\varepsilon) \sum^{q_1}_{i=0} \sum^{q_2}_{j=0} \beta^2_{ij}\, , \\
    &\pgf_X(u) = \prod^{q_1}_{i=0} \prod^{q_2}_{j=0} \pgf_\varepsilon\bbrackets{1+\beta_{ij}(u-1)}\,,
\end{split}
\end{equation}
where $\pgf_\varepsilon$ denotes the innovations' pgf. Note that similar to the CINAR model \eqref{eq: CINAR recursion}, but different from the ordinary INAR models for random fields, Poisson innovations go along with Poisson observations. 

\smallskip
Compared to the (C)INAR random fields, further independence assumptions on the ``past'' observations are not necessary in order to derive the model properties. In particular, dependence between those thinnings being applied to the \emph{same}~$\varepsilon_{s,t}$ is not prohibited (similar to the INMA time-series model). This means that conditioned on $\varepsilon_{s,t}$, the $(q_1+1)(q_2+1)$ random variables 
\[
    \beta_{00}  \circ_{s,t} \varepsilon_{s,t},\ \beta_{01}  \circ_{s,t+1} \varepsilon_{s,t},\ \dots,\ \beta_{q_1 q_2} \circ_{s+q_1,t+q_2} \varepsilon_{s,t}
\]
might exhibit cross-dependence, which leads to different model specifications that have different ACvFs but the same marginal properties. Let $(Z_{s,t;r}^{(i,j)})_{1\leq r \leq \varepsilon_{s,t}}$ denote the counting series involved in the thinning applied to $\varepsilon_{s,t}$ at point $(s+i,t+j)$ for $0 \leq i \leq q_1$ and $0 \leq j \leq q_2$, so $\beta_{ij}\circ_{s+i,t+j} \varepsilon_{s,t} = \sum_{r=1}^{\varepsilon_{s,t}} Z_{s,t;r}^{(i,j)}$ with $\bbp(Z_{s,t;r}^{(i,j)}=1)=\beta_{ij}$. Accordingly, the $(q_1+1)(q_2+1)$-dimensional vectors $\bfZ_{s,t;r}:= (Z_{s,t;r}^{(0,0)}, \dots, Z_{s,t;r}^{(q_1,q_2)})^\top$ are i.i.d.\ for all $s,t$ and $r$, but their components might be cross-dependent.
Defining the index sets $\mathcal{S}_{kl} := \{(i+k,j+l) : 0 \leq i \leq q_1, 0 \leq j \leq q_2\}$ for $k,l\in\bbz$, the ACvF satisfies 
\begin{align}
    \gamma(k,l) &= (\sigma^2_\varepsilon - \mu_\varepsilon) \sum_{(i,j)\in\mathcal{S}_{00}\cap\mathcal{S}_{kl}} \beta_{i,j}\beta_{i-k,j-l} \nonumber\\
    &\qquad + \mu_\varepsilon \sum_{(i,j)\in\mathcal{S}_{00}\cap\mathcal{S}_{kl}} \bbp(Z_{s-i,t-j;1}^{(i,j)}=Z_{s-i,t-j;1}^{(i-k,j-l)}=1) \label{eq: inma sacf}
\end{align}
(with the convention that empty sums equal $0$),
which covers the variance formula from \eqref{eq: INMA mean var pgf} for $k=l=0$. 

\smallskip
Inspired by the models known from time series analysis, \citet{sil_wei_26b} proposed two special cases concerning the dependence structure within $\bfZ_{s,t;r}$, which are well-interpretable in application: the \emph{independence model} assumes independence between the components of $\bfZ_{s,t;r}$ (hence, $\bbp(Z_{s-i,t-j;1}^{(i,j)}=Z_{s-i,t-j;1}^{(i-k,j-l)}=1) = \beta_{i,j}\beta_{i-k,j-l}$ for $(k,l)\not=(0,0)$ in \eqref{eq: inma sacf}), while the \emph{spread model} has $\bfZ_{s,t;r} \sim \text{MULT}(1; \beta_{00}, \dots, \beta_{q_1 q_2})$ with $\beta_\bullet\leq 1$, i.e., at most one component of $\bfZ_{s,t;r}$ is equal to one and all others are equal to zero (hence, $\bbp(Z_{s-i,t-j;1}^{(i,j)}=Z_{s-i,t-j;1}^{(i-k,j-l)}=1) = 0$ for $(k,l)\not=(0,0)$ in \eqref{eq: inma sacf}). 
Note that the Poisson INMA$(1, 1)$--independence model was applied to the yeast count data from Figure~\ref{fig: examples}\,(a) as a possible alternative to the INAR$(1,1)$ model from Sect.~\ref{section: The Origins: The First-order INAR Random Field}. 
Moreover, \citet{sil_wei_26b} note that \eqref{eq: model inma} naturally extends to a \emph{multilateral INMA$(q_1,q_2,p_1,p_2)$ model} like in Sect.~\ref{section: INAR Extensions to Multilateral Structure}, namely by considering the recursion
\begin{equation}
    \label{eq: model inma multilateral}
    X_{s,t} = \sum^{q_1}_{i=-p_1} \sum^{q_2}_{j=-p_2} \beta_{ij} \circ_{s,t} \varepsilon_{s-i,t-j}\, .
\end{equation}
They show that \eqref{eq: INMA mean var pgf} and \eqref{eq: inma sacf} still hold true after slight adaptions.

\section{Conclusions and Outlook}
\label{section: Conclusions and Outlook}

This article showed that during the last 15 years, an impressive variety of INARMA models has been developed for count random fields, covering approaches with different thinning operators, first- and higher-order models, as well as unilateral and multilateral model structures. These models are characterized by making use of thinning operators (most often binomial thinning) instead of ordinary multiplications, and they exhibit a linear spatial autocorrelation structure similar to that of classical ARMA models. But it also became clear that there are still several notable gaps in the research literature which, in turn, offer attractive opportunities for future research. For example, all model proposals refer to unbounded counts (i.e., with full~$\bbn_0$ as their range), whereas ARMA-like models for bounded counts (i.e., with range $\{0,\ldots,N\}$ for some $N\in\bbn$) have not been developed so far. An analogous gap can be identified with respect to signed integer outcomes (i.e., if~$X_{s,t}$ has range~$\bbz$). Also the possibility of negative dependencies deserves further research activity, e.g., by exploring the aforementioned Tobit approach in further detail. In addition, different ARMA-like models for count or integer time series, such as the INGARCH models (recall Section~\ref{section: Modifications Based on Alternative Thinning Concepts and Mixtures}) or rounded ARMA models, have not been adapted to the case of random fields so far. Finally, there are several open questions even within the existing model proposals, for example with respect to the multilateral INAR model, also recall Remark~\ref{remark:multilateral_open_questions}. Another interesting question refers to the ``orientation'' of the unilateral INAR models (e.g., if comparing $(X_{s,t})$ to $(X_{-s,t})$). For INAR time series models, time reversibility has been proven for a few instances, see \citet{schweer15}, but we are not aware of corresponding extensions to INAR random fields.

\subsection*{Acknowledgement}
The authors are grateful to the reviewer for useful comments on an earlier draft of this article.

\bibliographystyle{plainnat}
\bibliography{references}
\end{document}